\documentclass[dvips,twoside,dvips]{article}
\usepackage[a5paper,layoutwidth=148.5mm,layoutheight=7.77817in,
outer=10mm,inner=38pt,tmargin=10mm,bmargin=10mm,footskip=15pt,
]{geometry}

\usepackage[utf8]{inputenc}
\usepackage[LGR,T1]{fontenc}
\usepackage{mathptmx}

\usepackage[greek,french]{babel}
\def\grec#1{\foreignlanguage{greek}{#1}}

\usepackage{graphicx}

\usepackage[leftmargin=5mm,rightmargin=5mm]{quoting}

\usepackage[hyphens]{xurl}%

\usepackage[comma]{natbib}
\setcitestyle{aysep=}
\usepackage{doi}

\begin{document}
\title{Expériences de l'infini et de la généralité dans le premier livre des \emph{Éléments} d'Euclide}
\author{Stefan Neuwirth}
\date{}
\maketitle


Cette contribution vise à proposer une lecture du premier livre des \emph{Éléments} d'Euclide qui ressuscite quelques potentialités du texte en imaginant les choix mathématiques et discursifs qui ont pu le commander. Au préalable, je vais contextualiser cette lecture en précisant comment j'approche cette œuvre.

Quel but puis-je lui assigner? La présentation de la géométrie qu'elle expose a des prérequis: elle s'adresse à une personne qui a déjà pratiqué cette science, qui a déjà mené des lignes droites et décrit des cercles, qui envisage déjà la géométrie comme une science de l'incidence, de l'intersection et de la position relative (l'inclinaison) de lignes, de l'égalité par ajustement. Les lecteurs et auditeurs des \emph{Éléments} ont déjà pratiqué les diagrammes comme matérialisations des figures, avec des lettres pour référer à une figure ainsi qu'aux choses qu'elle comporte (les points, les lignes, les angles, les surfaces, voire les nombres dans les livres arithmétiques).

Pour saisir comment cette œuvre s'articule avec les connaissances et pratiques de ses lecteurs, je conçois que celles-ci sont en constante évolution et suivent leur propre chemin: elles donnent lieu à une grande variété d'approches et de rapports aux choses; si la découverte des \emph{Éléments} peut les affecter, il n'y a là cependant aucune nécessité. Cette articulation est pour moi le lieu de la pensée diagrammatique.

J'interprète le silence des \emph{Éléments} à l'égard des prérequis comme une licence à poursuivre son chemin et à développer son propre point de vue, comme une manière de tenir ouverte la possibilité de leur évolution. Les rédacteurs des \emph{Éléments} revendiquent néanmoins, par le soin des références internes au texte et par son architecture, l'autonomie de leur rédaction: c'est-à-dire qu'ils entendent proposer une œuvre qui est cohérente et valable en tant que telle, indépendamment des points de vue et de leur évolution. Proposer une lecture vivante, c'est pour moi d'abord essayer d'être à la hauteur de cette revendication en témoignant comment elle m'affecte tant dans ses silences que dans le flot de son discours.

Cette autonomie vaut en particulier vis-à-vis d'éventuelles problématiques de fondation. En ce qui me concerne, les \emph{Éléments} tiennent un discours qui vaut en tant que tel et qui résonne avec notre connaissance et notre pratique: celles-ci le précèdent chronologiquement dans la mesure où elles nous permettent de comprendre son intention; et cependant il ne vise pas le monde matériel ni les diagrammes, il vise au contraire sa propre forme et sa propre organisation.

Il y a donc un aspect prométhéen dans les \emph{Éléments}: par leurs silences et par l'abondance des mots, ils montrent qu'il suffit de parler à leur manière pour capter et former le savoir géométrique. Plutôt que les définitions, demandes, notions communes elles-mêmes, c'est la mise en place de ce discours visant à la généralité qui a imprimé sa marque aux sciences mathématiques.

Dans cette organisation, le langage formulaire, voire répétitif et incantatoire, le vocabulaire restreint font pendant à la richesse des formes grammaticales, à leur exploitation voire invention dans le cadre de la langue grecque puis des langues dans lesquelles les \emph{Éléments} ont été traduits. C'est aussi ce processus vivant de création dans la langue naturelle que je vise comme circonstance de leur invention. La langue elle-même, en résonance avec les connaissances et pratiques antérieures, devient ainsi le laboratoire d'Euclide.

\begin{sloppypar}
Ce discours invite à imaginer des dialogues qui amplifient sa signification: par exemple, la voix qui annonce « Je dis que\dots » comme un défi dans chaque démonstration s'adresse au lecteur qui veut être convaincu, qui est à la fois coopératif et cependant méfiant. Au cours de la démonstration qui suit, le langage formulaire permet alors de repérer les voix des assertions invoquées afin qu'elles apportent leur contribution. Réciproquement, notre connaissance de la forme du dialogue issue de la vie courante est investie dans notre lecture.
\end{sloppypar}

Euclide ne propose pas de réflexion mathématique sur la logique en général ni sur la quantification en particulier. Et pourtant le langage est à l'œuvre, d'abord pour assurer la généralité du discours, et il me semble utile de constater ceci et de le nommer à l'aide des concepts que je connais même s'ils sont postérieurs à la genèse des \emph{Éléments}.

Dans la suite, je vais me concentrer sur le premier livre de ce traité et témoigner de la variété de mes expériences de l'infini et de la généralité à sa lecture. Ce focus éclaire des facettes souvent négligées d'un sujet central des mathématiques grecques par l'exhaustion d'un corpus limité et ce dans l'étendue limitée du texte et dans le temps limité de la recherche impartis. Je considère cette étude comme un travail préliminaire à une recherche sur la cognition incarnée de l'infini \citep{lakoffnunez00,nunez15} et de ses ancres matérielles \citep{hutchins05}.

\bigskip

La lecture des \emph{Éléments} procure des expériences de l'infini pour des raisons qui relèvent soit de la géométrie elle-même soit du discours élaboré.

\section*{L'infini géométrique}

Je retrouve dans les \emph{Éléments} des choses que je comprends à l'aide de mon expérience du monde et que je rapporte à l'infini.

Partons d'une expérience du monde sensible. Lorsque nous appréhendons une chose, nous pouvons distinguer deux modalités.
\begin{itemize}
\item La chose s'offre à notre investigation dans son entier, comme un tout. Une partie de cette chose la délimite et en forme la frontière, en deçà de laquelle se trouve l'intérieur de la chose. L'intérieur s'offre de manière plus obscure, indéfinie, vague à l'investigation que la frontière.
\item La chose ne s'offre à nous que partiellement. Cependant, nous pouvons savoir comment poursuivre notre investigation de proche en proche: la partie connue donne accès à la partie encore inconnue et ainsi la chose est connaissable en puissance.
\end{itemize}
Dans la géométrie des \emph{Éléments}, cette distinction reste pertinente.

\subsection*{Les limitants et les illimités}

L'opposition des limitants aux illimités est la mieux attestée parmi celles qu'on attribue aux pythagoriciens par le rôle que lui donne Philolaos.
\begin{quoting}
  \noindent Ce sont les illimités et les limitants qui ont, en s'harmonisant,
  constitué au sein du monde la nature, ainsi que la totalité du monde
  et tout ce qu'il contient \citep[fragment B\,\textsc{i},][p.~502]{dumont88}.
\end{quoting}

\subsubsection*{La limite et l'intérieur}

Les définitions 13~et~14 disent ceci \citep[je vais utiliser dans tout l'article la traduction de Bernard][]{vitrac90}.
\begin{quoting}
  \noindent Une \emph{frontière} est ce qui est limite de quelque chose.
    
  \noindent Une \emph{figure} est ce qui est contenu par quelque ou quelques frontière(s).
\end{quoting}
Auparavant, les définitions 6~et~3 disent cela.
\begin{quoting}
  \noindent Les \emph{limites d'une surface} sont des lignes.
  
  \noindent Les \emph{limites d'une ligne} sont des points.
\end{quoting}
Vincenzo \citet{derisi} propose de donner une plus grande ampleur à la notion de limite: il lit aussi dans cette dernière définition l'expression euclidienne de l'existence de points d'intersection de lignes comme ceux des deux cercles considérés dans la proposition~1.

Puis les définitions 15, 18 et~19 définissent des figures particulières (cercle, demi-cercle, figure rectiligne) par leur frontière.
\begin{quoting}
  \noindent Un \emph{cercle} est une figure plane contenue par une ligne unique par rapport à laquelle [\dots].
    
  \noindent Un \emph{demi-cercle} est la figure contenue par le diamètre et la circonférence découpée par lui [\dots].
    
  \noindent Les \emph{figures rectilignes} sont les figures contenues par des droites [\dots].
\end{quoting}

Pour désigner un point de l'intérieur d'une chose, par exemple sur le prolongement d'une ligne droite, Euclide écrit ceci dans la proposition~5.
\begin{quoting}
  \noindent[\dots] qu'un point F soit pris au hasard [\grec{τυχὸν}, « comme il vient »] sur BD [\dots].
\end{quoting}
Pour désigner un point dans un demi-plan au delà d'une droite, il écrit cela.
\begin{quoting}
  \noindent[\dots] que soit pris au hasard le point D, de l'autre côté de la droite AB [\dots].
\end{quoting}

Le participe \grec{τυχὸν} maintient aussi le dialogue avec le lecteur parce que celui-ci peut se sentir appelé à placer un point arbitrairement et à mettre ce faisant l'auteur au défi de poursuivre la démonstration.

Si je place cette observation ici alors que les choses considérées sont illimitées, c'est que j'attribue l'usage de l'expression « comme il vient » non pas au fait que la droite et le demi-plan en question peuvent être prolongés à l'infini, mais au fait que l'intérieur d'une figure est illimité au sens qu'un point n'y a pas un rôle différent d'un autre: les points de l'intérieur sont génériques. Mon attribution est confirmée par les propositions du deuxième livre.

\subsubsection*{Les illimités}

La demande~3 dit ceci.
\begin{quoting}
  \noindent[Qu'il soit demandé] de prolonger continument en ligne droite une ligne droite limitée.
\end{quoting}
Si dans une forêt je veux avancer en ligne droite devant moi, je vais viser la direction de deux arbres qui se superposent dans mon regard. Si je veux prolonger la ligne d'un fil tendu entre deux piquets au delà de l'un d'eux, je vais attacher une extrémité d'un fil à celui-ci et le tendre en plantant un troisième piquet dont je vais viser l'emplacement en alignement avec les deux premiers piquets. Je peux vérifier mon choix en tendant un autre fil entre les piquets extrêmes pour constater qu'il se superpose tant au fil prolongé qu'au fil qui prolonge.

La demande~3 garantit la possibilité d'une telle démarche en géométrie. La ligne droite non limitée est donnée comme une amorce et c'est la demande~3 qui garantit la poursuite de la droite au delà de l'amorce. Dans ce sens, la ligne droite illimitée est une chose qui, de l'un ou des deux côtés, est donnée avec l'intention d'être prolongée: elle est donc donnée comme infinie au sens où on peut encore la prolonger exactement de la manière dont les nombres entiers sont infinis. Le lecteur qui lit la conclusion de la demande~5,
\begin{quoting}
  \noindent les deux droites, indéfiniment prolongées, se rencontrent du côté où sont les angles plus petits que deux droits,
\end{quoting}
peut la comprendre ainsi: si on prolonge les deux droites suffisamment loin, elles se couperont, c'est-à-dire que l'existence d'une intersection résulte de ces prolongements. Cette lecture d'un infini en puissance me comble.

En conclusion, dans l'opposition des choses limitées aux choses illimitées, je soutiens qu'Euclide réfère à la manière dont l'objet est donné et non pas à une propriété intrinsèque de l'objet.

\subsection*{L'infini par addition et l'infini par division en géométrie}

Ces deux infinis figurent parmi les infinis des mathématiciens dans le troisième livre de la \emph{Physique} d'Aristote.
\begin{quoting}
  \noindent La croyance qu'il existe un certain infini pourrait se tirer, pour qui examine la question, de cinq arguments principaux: [\dots]; de la division dans les grandeurs (en effet, les mathématiciens eux aussi ont recours à l'infini); [\dots]. Mais il y a surtout le principal argument qui offre à tous la même difficulté: en effet, du fait qu'on ne leur trouve pas de borne dans la pensée, on est d'avis que le nombre, aussi bien que les grandeurs mathématiques et ce qui est hors du ciel sont infinis \citep[\emph{Physique}~203\,b\,15,][p.~176]{pellegrin00}.
\end{quoting}
L'expérience de l'infini par addition se fait dans la respiration rythmée par l'expiration et l'inspiration, dans le pouls rythmé par les battements du cœur, dans l'itération du pas de la marche. L'expérience de l'infini par division se fait dans le partage: partage d'une boisson, d'une denrée, d'un butin.

Ces deux infinis s'incarnent dans les choses géométriques.

\subsubsection*{L'infini par addition}

\begin{figure}[h]
  \centering
  \includegraphics[scale=.8]{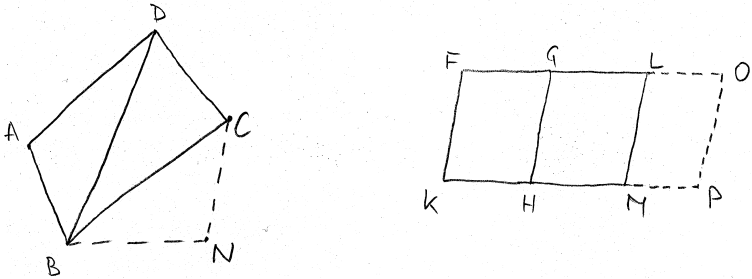}
  \caption{Construire un parallélogramme égal à une figure rectiligne donnée.}
  \label{121}
\end{figure}

La définition~19 dit ceci.
\begin{quoting}
  \noindent Les \emph{figures rectilignes} sont les figures contenues par des droites; \emph{trilatères}: celles qui sont contenues par trois droites, \emph{quadrilatères} par quatre; \emph{multilatères} par plus de quatre.
\end{quoting}
L'infini des nombres se retrouve ainsi dans le nombre de lignes droites qu'il peut falloir pour contenir une figure rectiligne. La proposition 45 énonce ceci.
\begin{quoting}
  \noindent Dans un angle rectiligne donné, construire un parallélogramme égal à une figure rectiligne donnée.
\end{quoting}
La figure donnée est donc contenue par un nombre donné de lignes droites. Or la démonstration est faite pour un quadrilatère, ABCD, qui est triangulé en 2~triangles, ABD, DBC (figure~\ref{121}): un parallélogramme égal au 1\ier\ triangle est construit dans l'angle donné par la proposition~42, FGHK; puis un 2\ieme\ parallélogramme égal au 2\ieme\ triangle est « appliqué » sur le 1\ier\ parallélogramme par la proposition~44, GHML, de sorte que les 2 ensemble forment encore un parallélogramme, FKML, qui est par conséquent égal au quadrilatère.
Euclide traite donc le cas particulier d'un quadrilatère contenu par 4~lignes droites. Le lecteur doit se rendre compte de cela tout seul; protester intérieurement; se poser la question de ce qui se passe lorsqu'on passe à une figure contenue par 5~lignes droites; se répondre qu'elle est triangulée en 3~triangles; appliquer 2~fois des parallélogrammes de sorte qu'on obtienne successivement des parallélogrammes qui sont égaux à~2 puis aux 3~triangles; c'est alors qu'il se rendra compte que ce qu'Euclide a dit, « le dire une fois revient à le dire sans cesse » \citep[Zénon d'Élée, fragment B\,\textsc{i},][p.~291]{dumont88} et que les applications successives de triangles ne diffèrent pas de celle qu'il a décrite, et que dans ce sens la démonstration du cas particulier est une démonstration du cas général.

Toute cette réflexion sur l'itération se fait face au silence d'Euclide et dans l'imagination du lecteur, qui supplée la nécessité de l'itération et reconnait le caractère paradigmatique du cas traité.

\subsubsection*{L'infini par division}

La définition~15 dit ceci.
\begin{quoting}
  \noindent Un \emph{cercle} est une figure plane contenue par une ligne unique par rapport à laquelle toutes les droites menées à sa rencontre à partir d'un unique point parmi ceux qui sont placés à l'intérieur de la figure, sont égales entre elles.
\end{quoting}
Elle dit donc que les rayons d'un cercle sont tous égaux entre eux. Or ceux-ci forment un infini par division: en effet, entre deux rayons il y a encore un rayon, construit par exemple à l'aide de la proposition~9.

\subsubsection*{Le mouvement}

Mais l'infini de la physique le plus mystérieux, celui évoqué par les paradoxes de Zénon, tant ceux du mouvement que ceux de la pluralité \citep{neuwirth21}, est tu par Euclide. Ce silence se manifeste pour moi dans les demandes 1, 2~et~3.
\begin{quoting}
  \noindent Qu'il soit demandé de mener une ligne droite de tout point à tout point.

  \noindent Et de prolonger continument en ligne droite une ligne droite limitée.

  \noindent Et de décrire un cercle à partir de tout centre et au moyen de tout intervalle.
\end{quoting}
Il se manifeste encore dans la notion commune~7.
\begin{quoting}
  \noindent Et les choses qui s'ajustent les unes sur les autres sont égales entre elles.
\end{quoting}
Ces demandes expriment que la géométrie d'Euclide commence là où les paradoxes s'arrêtent: elles mettent en suspens la perplexité qu'ils causent en demandant de concéder ces constructions
. Elles présentent les actes de mener, de prolonger et de décrire comme des oracles, comme des apports extérieurs à la géométrie: lorsqu'une figure propose deux points, l'invocation de la demande~1 fera mener une ligne droite d'un point à l'autre d'une manière qui reste inexprimée. C'est au lecteur de suppléer son expérience: par exemple,
\begin{itemize}
\item planter un piquet puis un autre en terre, attacher une
  extrémité d'un  fil au premier piquet, se rendre au deuxième piquet, tendre le fil et l'y attacher, constater la présence de ce fil tendu;
\item planter un piquet en terre, décider d'un intervalle de fil, attacher une
  extrémité du fil au piquet, s'en éloigner pour tendre le fil, marquer la ligne décrite par l'autre extrémité lorsqu'on se déplace d'une manière qui laisse le fil tendu en se rendant compte qu'il y a deux directions possibles et que pour chacune on revient à la position initiale au bout d'un certain temps, constater la présence de cette marque.
\end{itemize}
Ces expériences de l'étendue sont contraintes par les circonstances terrestres: les lignes droites menées à la surface de la Terre sont courtes; la ligne d'Alexandrie à Syène considérée par Ératosthène ne peut pas être décrite à partir du centre de la Terre.

\begin{figure}[h]
  \centering
  \includegraphics[scale=.699]{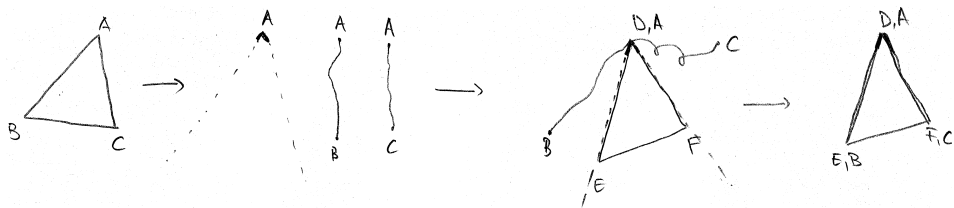}
  \caption{Ajuster un triangle sur un autre.}
  \label{123}
\end{figure}

J'associe l'expérience suivante à la notion commune~7:
\begin{itemize}
\item mesurer un fil par un autre: pincer ensemble une des extrémités de chaque fil entre le pouce et l'index de la main gauche; appliquer un fil contre l'autre au voisinage immédiat de cette extrémité entre le pouce et l'index de la main droite; tirer les fils à partir de leur extrémité tenue par la main gauche tout en les laissant glisser entre les doigts de la main droite de telle sorte que pendant tout le processus les fils restent appliqués l'un contre l'autre en ce lieu au sens qu'ils n'y glissent pas l'un relativement à l'autre; atteindre ce faisant l'autre extrémité de l'un des deux fils; constater qu'ainsi ce fil a été ajusté sur une partie de l'autre, qui lui est donc égale.
\end{itemize}
Paraphrasons l'ajustement d'un triangle sur un autre, ABC sur~DEF, dans la proposition~4 (figure~\ref{123}): on pose l'angle sous~A sur l'angle sous~EDF; on constate que le point~A est posé sur~D et que les directions des angles coïncident, chacune à chacune; on tire le côté~AB le long du côté~DE et on constate que ce faisant le point~B est posé sur~E; on tire le côté~AC le long du côté~DF et on constate que ce faisant le point~C est posé sur~F. 

Les expériences que j'ai proposées relèvent de l'étendue et du temps, que nous con\-nais\-sons sans pour autant savoir comment ils sont faits: la géométrie des \emph{Éléments} procède de manière à ne pas éprouver cette ignorance.

\section*{La généralité du discours}

Le lecteur des \emph{Éléments} rencontre tout au long du traité le discours qui assure la généralité. Voici un aperçu de la variété linguistique dans laquelle il se manifeste dans la langue grecque. Il s'appuie sur les analyses par Fabio \citet[p.~122]{acerbi11} de ce qu'il nomme les « déterminants de généralité »:
\begin{itemize}
\item l'indétermination:
  \begin{itemize}
  \item l'absence d'article (rendue si besoin par l'article indéfini en français): la définition~1 dit « Un \emph{point} est ce dont il n'y a aucune partie »;
  \item l'adjectif \grec{τις} (traduit par « n'importe quel » par Bernard \citealt{vitrac90} et par « un certo [quelque] » par Fabio \citealt{acerbi07}): la définition~17 dit « Et un \emph{diamètre} du cercle est n'importe quelle droite menée par le centre, limitée de chaque côté par la circonférence du cercle, laquelle coupe le cercle en deux parties égales »;
  \end{itemize}
\item la totalité (quantification universelle):
  \begin{itemize}
  \item l'adjectif \grec{πᾶς} [tout]: la définition~15 du cercle dit que « toutes les droites menées à sa rencontre à partir d'un unique point parmi ceux qui sont placés à l'intérieur de la figure, sont égales entre elles » et il vient avec un potentiel infini d'égalités qui se déploie à mesure que l'on introduit ses intersections avec d'autres lignes; la demande~1 est « de mener une ligne droite de tout point à tout point »; la proposition~16 énonce que « dans tout triangle, un des côtés étant prolongé, l'angle extérieur est plus grand que chacun des angles intérieurs et opposés »;
  \item les articles déterminatifs \grec{τῆς} et \grec{ὸ} (rendus par l'article indéfini): la proposition~1 énonce « sur une [\grec{τῆς}, « la »] droite limitée donnée, construire un triangle équilatéral »; Fabio \citet[p.~154]{acerbi11} explique leur usage  à la fois comme un choix stylistique que par le contexte « dialectique, où qui interroge propose un objet géométrique et qui répond a la tâche d'y effectuer la construction »;
  \item l'article déterminatif au pluriel: la proposition~5 énonce que « les angles à la base des triangles isocèles sont égaux entre eux » et réfère ainsi à des classes entières d'objets.
  \end{itemize}
\end{itemize}
La signification de la généralité se construit dans l'apprentissage de la langue naturelle, qui utilise les mêmes moyens d'expression et dans laquelle nous repérons d'emblée les assertions qui en relèvent et celles qui n'en relèvent pas. Je propose ci-dessous des exemples tirés du \emph{Dictionnaire} de l'Académie française (\citeyear{dictionnaire35}, \citeyear{dictionnaire92}) en suivant l'ordre de l'aperçu ci-dessus.
\begin{itemize}
\item « Je cherche un livre sur tel sujet. »
\item « Prenez n’importe quel livre dans ma bibliothèque. »\\
  « Si cela était, quelque historien en aurait parlé. »
\item « Toute sa famille est en bonne santé. »
\item « L’homme, la femme. » « Le vin de Bourgogne. » « Il vient le mardi. »
\item « Les livres que vous m’aviez prêtés. »
\end{itemize}

\subsection*{L'exposition des propositions}

« L'exposition isole le donné en lui-même et le dispose pour qu'il puisse permettre la recherche », écrit Proclus dans son \emph{Commentaire au premier livre des \emph{Éléments} d'Euclide} \citep[voir][p.~279]{guillaumin06}. Elle commence dans la grande majorité des cas de la même manière, traduite ainsi par Vitrac pour les propositions~16, 17, 39~et~40.
\begin{quoting}
  \noindent Soit le triangle~ABC [\dots].
  
  \noindent Soit un triangle~ABC [\dots].
  
  \noindent Soient ABC,~DBC des triangles égaux [\dots].
  
  \noindent Soient les triangles égaux ABC,~CDE [\dots].
\end{quoting}
Il y a cependant des variations. Le « donné » peut être le sujet d'un verbe: faire, produire, se couper, tomber, avoir la même base, être dans les mêmes parallèles, être égal; il peut être renvoyé dans une subordonnée.

La généralité mathématique ainsi installée dans les démonstrations des \emph{Éléments} d'Euclide est l'objet de discussions par de nombreux auteurs, comme le relève Fabio \citet[note~2, p.~28]{acerbi20}. Celui-ci critique les traductions ci-dessus, qu'il amende ainsi.
\begin{quoting}
  \noindent Soit un triangle, ABC [\dots].

  \noindent Soient des triangles égaux, ABC,~DBC [\dots].
\end{quoting}
En effet, il note que le verbe être a ici une valeur non pas existentielle mais présentielle, c'est-à-dire que sa seule fonction est d'introduire la liste des choses données dans les hypothèses de la proposition. De cette manière, la généralité de l'énoncé et celle de l'exposition coïncident: l'un et l'autre introduisent simplement les choses dont la proposition va traiter.

\subsection*{La quantification dans la définition et le postulat des parallèles}

La définition des droites parallèles est la seule à avoir une forme négative, la non-existence d'un point d'intersection de ces droites. C'est pourquoi, dans la proposition~27, le parallélisme est démontré en déduisant une contradiction de l'existence d'un point d'intersection des droites. Alors que la non-existence nécessite l'évocation d'un infini, le prolongement indéfini des droites considérées, l'existence, elle, correspond à une situation finie: deux droites avec leur point d'intersection.

\begin{figure}[h]
  \centering
  \includegraphics[scale=.8]{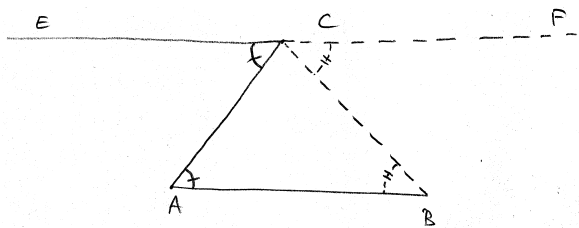}
  \caption{La somme des angles d'un triangle.}\label{212}
\end{figure}

La demande~5 (le postulat des parallèles) postule un critère pour l'existence d'un point d'intersection de deux droites et donc pour leur non-parallélisme. La proposition~29 formule et établit la contraposée du postulat. Celle-ci illustre que la demande~5 réduit une infinité d'inclinaisons possibles de deux droites parallèles à une seule possibilité: qu'une droite quelconque tombant sur des droites parallèles fasse « l'angle extérieur égal à l'angle intérieur et opposé », ni plus ni moins. C'est cela qui permet d'établir la proposition~32 sur la somme des angles d'un triangle (figure~\ref{212}): dans sa démonstration, la parallèle à~AB passant par~C est construite selon la proposition~31 par rapport à la droite~AC, c'est-à-dire qu'un point~E est introduit de sorte que l'angle sous~ACE soit égal à l'angle sous~CAB; la proposition~29 montre que les points~C et~E sont en alignement avec le point~F qui serait introduit pour construire la parallèle par rapport à la droite~BC, c'est-à-dire de sorte que l'angle sous~BCF fût égal à l'angle sous~CBA.

\subsection*{La quantification dans les raisonnements par l'absurde}

Dans le premier livre des \emph{Éléments}, les raisonnements par l'absurde visent chacun à démontrer une égalité (d'angles) ou une coïncidence (de points). Pour démontrer ces égalités ou coïncidences, Euclide va démontrer l'absurdité de la différence en étudiant l'hypothèse que deux grandeurs ou choses sont différentes; or la différence recèle une infinité de possibilités. Cela appelle deux commentaires. 
\begin{enumerate}
\item Comme le note Henri \citet{lombardi97}, l'égalité propre aux rapports de grandeurs est la négation de leur inégalité. Cet état de fait se retrouve dans les définitions~7 et 5 du cinquième livre: l'inégalité de deux rapports de grandeurs résulte de l'existence de multiples des grandeurs avec une certaine propriété, alors que l'égalité requiert la considération de tous les multiples. 
\item Dans les propositions~5, 6 et~14, Euclide pourrait constater une égalité, mais il ne le fait pas; elle proviendrait de la coïncidence d'une chose avec elle-même; lorsqu'il constate une égalité, c'est pour des choses distinctes.
\end{enumerate}

Dans la proposition~7, cette infinité est introduite par la considération de « points différents, C~et~D [\grec{ἄλλῳ σημείῳ τῷ τε Γ καὶ Δ}] » (la traduction de Vitrac, « un point quelconque~D, différent de~C », casse la symétrie de l'introduction des deux points).

\begin{sloppypar}
Dans la proposition~14, la désignation de cette infinité a lieu d'une autre manière: pour montrer que BD est en alignement avec BC, Euclide va supposer qu'elles ne le sont pas, et alors introduire la droite BE en alignement avec BC pour déduire une contradiction. L'argument se termine ainsi.
\begin{quoting}
  \noindent Donc BE n'est pas en alignement avec CB. Alors semblablement nous démontrerons qu'une autre droite quelconque [\grec{ἄλλη τις}, « quelque autre »] ne l'est pas non plus sauf BD. Donc CB est en alignement avec BD.
\end{quoting}
L'infini apparait donc tardivement, dans la considération d'« une autre droite quelconque ». Cette infinité d'autres droites est alors contredite par la généralité de la démonstration proposée, qui permet de concevoir des démonstrations « virtuelles » (\citealp[note~85, p.~224]{vitrac90}, appelées « potentielles » par \citealp[\S~1.5.2]{acerbi11})  pour chacune d'elles.
\end{sloppypar}


\end{document}